\documentclass [11pt]{article}
\usepackage{amstext}

\input{artadj.sty}
\usepackage{mathtools}   %bmatrix
\usepackage{amsmath}
\usepackage{enumitem}  %enumerate list
\usepackage{latexsym}  
\usepackage{graphicx}
\usepackage{longtable}
\usepackage{listings}
\usepackage{setspace}

\renewcommand{\thesection}{\arabic{section}}

\newtheorem{thm}{Theorem}[section]
\newtheorem{lem}[thm]{Lemma}

\newcommand {\done} {\quad\vrule height4pt WIDTH4PT}

\newcommand{\inof}{\int_{0}^{\infty}}
\newcommand{\bos}{\boldsymbol}
\begin{document}

\title{
An Efficient Method to Compute the Stationary Probabilities of the $GI^X/M/c/N$ Model}
\author{ 
Muhammad El-Taha\\ 
email:el-taha@maine.edu\\
and\\
Thomas Michaud\\
email:thomas.michaud@maine.edu\\
\\
Department of Mathematics and Statistics\\
University of Southern Maine\\
96 Falmouth Street\\
Portland, ME  04104-9300\\
}
\date{ }

\maketitle

\ls{1}

\abstract{
Consider the batch-arrival $GI^X/M/c/N$ model with $c$ servers, general inter-arrival batch times, finite buffer, and exponential service times. Inter-arrival batch times, batch sizes,  and service times are $i.i.d.$ and independent of each other. In this article we give a simple efficient algorithm to derive an exact solution for the steady state system size probabilities. The starting point is computing the one-step transition probabilities of the imbedded Markov chain observed at the system arrival epochs of the the corresponding  $G/M/c$ model. The one-step transition probabilities are computed exactly by converting a numerical integration problem into a finite sum.  Another key contribution is generating the transition probabilities of the batch-arrival model by using a simple and intuitive method to extend the results of the standard $GI/M/c$ model to batch arrivals with and without a finite buffer, and in the case of finite buffer with partial and full batch rejection. Moreover, we develop an efficient stable algorithm  that can  accurately solve problems with a larger number of servers  than  previously known.   We give numerical examples to demonstrate the performance of our method.}  
\bigskip

\noindent


\ls{1.3}
%%%%%%%%%%%%%%%%%%%%%%%%%%%%%%%%%%%%%%%%%%%%%%%%%%%%%%%%%%%%%%%%%%%%%%%%%%%%%%%%%%%%%%%%%%%%%%%%%%%%%

\section{Introduction}

In this article we study the performance characteristics of the  $GI^X/M/c/N$ multi-server queueing model with batch arrivals, general batch inter-arrival times, multiple parallel exponential servers, and finite buffer. We accomplish this in two steps. First, using the $GI/M/c$ model, we develop a simple efficient algorithm to accurately compute the one-step transition probabilities of the imbedded Markov chain where the imbedding points are the arrival epochs. We solve the imbedded Markov chain to compute the arrival-times stationary distribution and use that to compute the time-average stationary distribution. Second, we use an innovative  approach to extend the results of $GI/M/c$ queueing model to the $GI^X/M/c/N$ multi-server batch arrival finite capacity model. We cover the finite and the infinite buffer cases, and  for the finite buffer case we include models with partial and full batch rejection.  Specifically, in this article we give a remarkably  simple algorithm to derive an accurate solution for the steady state system size probabilities and give numerical examples to demonstrate the performance of our method.

%Literature on the $G/M/c$ and $G/M/c/N$
The well known multi-server  $GI/M/c$ queue is an important model in queueing theory.
This queueing model, with and without finite buffer, is studied by several authors.  Tak\'acs~\cite{Tak57} and ~\cite{Tak62} is, perhaps, the first author to study the $G/M/c$ model using the imbedded Markov chain approach. He develops a computational recursion that uses generating functions techniques. This method is described by Gross and Harris~\cite{Gro08} and
Kleinrock~\cite{Kle75} as "extremely long" and "complex". All textbooks that deal with this model  require the numerical integration of at least one transition probability expression of the imbedded Markov chain. See for example Gross and Harris~\cite{Gro08}, Kleinrock~\cite{Kle75}, Medhi~\cite{Med03}, Ross~\cite{Ros07}, and Tijms~\cite{Tij03}.  Others use approximation methods like Cosmetatos and  Godsave~\cite{Cos80}. Moreover, 
using the supplemental variable approach  Hokstad~\cite{Hok75} studies the $GI/M/c$ model with finite waiting room, and provides relations between  pre-arrival, post departure, and time-average probabilities. See also Yao et al~\cite{Yao84} who uses level crossing methods to relate  similar quantities.
Neuts~\cite{Neu81} uses a matrix analytic approach to investigate this and related models. Grassmann and Tavakoli \cite{Gra14} review  multiple  numerical approaches and address stability issues.	  
Kim and Chaudhry ~\cite{Kim17}  study the finite capacity  $GI/M/c/N$  model. See also Ferreira and Pacheco  \cite{Fer06}.

%Literature on Batch arrivals

The batch arrival  $GI^X/M/c$ model is studied by Chaudhry and Kim~\cite{Cha16}. They start with the balance equations for the imbedded Markov chain, write the characteristic equation, solve for  the  characteristic equation roots using  $MAPLE$, and use these roots to compute the arrival probabilities.
We note an oversight in  Chaudhry and Kim~\cite{Cha16} in the formula  cited at the bottom of page 241.  It is evident that the third and fourth expressions depend on $i$ and $j$ explicitly and not just the difference as indicated in the l.h.s.  $k_{i+h-j}$ of the equation.
A supplemental variable approach is used by  Laxmi and Gupta~\cite{Lax00}  to solve for the finite capacity $GI^X/M/c/N$ model in order to relate  pre-arrival  and time-average probabilities. They then use the imbedded Markov chain approach to obtain the pre-arrival probabilities, focusing only  on inter-arrival distribution functions  whose Laplace transform  can be analytically expressed like the Erlang and the hyper-exponential distribution functions.   See also Bailey and Neuts ~\cite {Bai81} who use matrix analytic methods for the study of the batch arrival $GI^X/M/c$ model. 
All these  methods need the one step transition probabilities stated in Lemma\ref{lem:pij}, but the  transition probabilities in  Lemma\ref{lem:pij} (iii) are approximated using numerical integration methods.
 %

%Contribution of this article: 
Computing the one step transition probabilities in the case of a transition from a state $i\geq c$ to a state $j\leq c-1$ requires tedious numerical integration. This is the key step in most approaches in the literature. In this article, the first contribution is a result that converts this expression into a finite sum. By itself, this result will make most approaches used in the literature function more efficiently. This and related results facilitate the development of a simple stable  algorithm to  efficiently compute the pre-arrival and time average probability distributions.  Another key contribution is extending  the $G/M/c$ standard multi-server model results to the $GI^X/M/c/N$ batch arrival model with and without finite capacity. We use a novel method to transform the transition probabilities of the standard  multi-server model to the batch arrival model. Furthermore, our approach in relating pre-arrival probabilities and time-average probabilities  uses intuitive and simple to follow sample-path methods.  Numerical results show that we can solve problems larger than previously known. 
We note that El-Taha~\cite{Elt21} uses a convolution approach to derive a result similar to Lemma~\ref{lem:pij}, however our current simpler direct method uses integration techniques only. Interestingly, our direct method complements the convolution approach given by El-Taha~\cite{Elt21}

Matrix analytic methods pioneered by Neuts~\cite{Neu81} are used to study this model and its variants. The solutions approach uses   algorithmic methods that utilize  matrix algebra and vectors. See also Bailey and Neuts ~\cite {Bai81} who use matrix analytic methods for the study of the batch arrival $GI^X/M/c$ model.

The rest of the article is organized as follows.  In Section~\ref{sec:sd} we focus on the standard $GI/M/c$ model and determine the transition probabilities, the system size probabilities at arrival instants, the time-average  stationary distribution,  and  measures of 
performance.  Specifically, we give  an efficient  algorithm to compute the pre-arrival and time-average stationary distribution.  In Section~\ref{sec:extensions} we give an intuitive method to extend our $GI/M/c$ results to the $GI^X/M/c$ and  $GI^X/M/c/N$ models. For the finite buffer model, we cover both the partial and full batch rejection.  In Section~\ref{sec:app} we give  numerical examples. The examples focus on the $GI/M/c$ model with and without finite buffer.  In the Appendix  we describe in detail the algorithm to compute the stationary probabilities.

\section{The $GI/M/c$ Model}\label{sec:sd}

In this section we focus on the  $GI/M/c$ model and observe the system at pre-arrival instants to determine the one step transition probabilities, and the system size probabilities, then the time average  stationary distribution is determined using the imbedded probability distribution. Measures of performance  can then be easily evaluated.  
The   $GI/M/c$ model  can be described as follows: We have $i.i.d.$ inter-arrival times $A_i, i\geq 1$ and $i.i.d.$ exponential service times $B_i, i\geq 1$  with common distribution functions $A(t)$ and $B(t)$ respectively. The first moments of the inter-arrival and service times are given by $E(A)=1/\lambda$ and $E(B)=1/\mu$ respectively. Note that for this model $\mu_n = \min(n,c)\mu$ is the state dependent service rate that shall be needed later.  The system state at pre-arrival epochs is described by the Markov chain $\{X_n, n=0,1,\cdots\}$ with transition probabilities given by  several standard textbooks in queueing theory.  Because these transition probabilities play an important role in our analysis and for ease of access we reproduce them  in the lemma below  (see for example Gross and Harris~\cite{Gro08} and El-Taha~\cite{Elt21}). 
\begin{lem}\label{lem:pij}
(i) for $ j\leq i+1\leq c$, the transition probabilities of $p(i,j)$ are given by 
\[
p(i,j)= \int_0^{\infty} \binom{i+1}{i-j+1} e^{-\mu tj}(1-e^{-\mu t})^{i-j+1}dA(t)\;; 
\]
(ii) For $c\leq j\leq i+1$, the transition probabilities of $p(i,j)$ are given by 
\[
p(i,j)= \int_0^{\infty} \frac {e^{-c\mu t}(c\mu t)^{i-j+1}}{(i-j+1)!}dA(t);
\]
(iii) for  $j+1 \leq c\leq i$
\[
p(i,j)= \binom{c}{c-j}\frac{(c\mu)^{i-c+1}}{(i-c)!}\int_0^{\infty}\int_0^{t}
  v^{i-c} e^{-\mu (t-v)j-c\mu v}(1-e^{-\mu (t-v)})^{c-j}dvdA(t)\;; 
\]
and $p(i,j)= 0$, otherwise.
\end{lem}

%Note that, in Lemma~\ref{lem:pij} (ii), if we let $b_n= \int_0^{\infty} \frac {e^{-c\mu t}(c\mu t)^{n}}{n!}dA(t),$  then we may write
%$p(i,j)=b_{i-j+1}\;.$

In the next subsection we provide our first key result. Specifically, we use direct integration to show that the one-step transition probability in Lemma~\ref{lem:pij} (iii) can be written as finite sum. This result complements a convolution method used by El-Taha~\cite{Elt21}.

\subsection{The One Step Transition Probabilities}

We start by defining the Markov chain imbedded at the pre-arrival instants, then describing our method for computing the transition probabilities.
Let $X_n$ be the number of customers  in the system at pre-arrival time instants, and $D_n$ be the number of customers served  during the $n^{th}$ inter-arrival time. 
Then, $X_{n+1}$ and $X_{n}$ are related by
\[
X_{n+1}= \left\{ \begin{array}{cc} 
                X_{n}+ 1 -D_n & \mbox{ $ D_n \leq X_{n}+1, X_n \geq 0$}\\
                0             & otherwise
                \end{array}
                \right. .\\
\] 
%where $D_n$ is the number of customers served  during the $n^{th}$ inter-arrival time. 
%
It is  straightforward to see that $\{X_n, n \geq 1\}$  is a Markov chain with one step transition probabilities  defined as 
$p(i,j)=P\{X_n =j|X_{n-1}=i\}, i=0, \ldots,; j=0, \ldots,$ and given by Lemma~\ref{lem:pij}. 
Let the Laplace-Stieltjes transform (LST) of the inter-arrival times distribution function $A(t)$ be given by $A^*(s)=\inof e^{-st}dA(t)$. Moreover, let $\frac{d^nA^*(s)}{ds^n}$ be the $n^{th}$ derivative  of $A^*(s)$ and denote $A^*_n(s)= (-1)^n \frac{d^nA^*(s)}{ds^n}$.
It can be easily verified that 
$A^*_n(s)=\inof t^n e^{-st}dA(t)$ for all $n\geq 0$, where $A^*_0(s)=A^*(s)$.
Now we present our  fundamental  result.
\begin{lem}\label{lem:fund} 
For $j >0, j+1 \leq c\leq i$ the one-step transition probabilities are given by
\begin{align} \label{eq:direct}
p(i,j) =\sum_{k=1}^{c-j} \frac{(-1)^{c-j-k}(c-1)!}{(k-1)!(c-j-k)!j!}  \left(\frac{c}{k}\right)^{i-c+2} \bigg[A^*((c-k)\mu ) - \sum_{r=0}^{i-c+1} \frac{(k\mu)^r A^*_r(c\mu)}{r!}  \bigg].
\end{align}
\end{lem}
{\bf Proof.}
To prove (\ref{eq:direct}) we  work with a modified version of  Lemma~\ref{lem:pij} (iii).
Define $\nu$ as the time required for $i-c+2$ service completions, and $c-j-1$ service completions for the remaining $t-\nu$ time. This leads to the equivalent equation
\begin{align*}
p(i,j) &= 
\int_{0}^{\infty}\!\!\!\! \int_{0}^{t} \binom{c-1}{c-j-1} e^{-\mu (t-\nu)j} (1-e^{-\mu(t-\nu)})^{c-j-1} \frac{(c\mu)^{i-c+2}\nu^{i-c+1}e^{-c\mu\nu}}{(i-c+1)!} \,d\nu\,dA(t).
\end{align*}
\begin{align*}
&= 
\binom{c-1}{c-j-1} \frac{(c\mu)^{i-c+2}}{(i-c+1)!} \int_{0}^{\infty}\!\! e^{-\mu jt} \int_{0}^{t}  e^{-\mu(c-j)\nu} \nu^{i-c+1}  (1-e^{-\mu(t-\nu)})^{c-j-1} \,d\nu\,dA(t)\\
\end{align*}
Note that  
\[(1-e^{-\mu(t-\nu)})^{c-j-1} = \sum_{k=0}^{c-j-1} \binom{c-j-1}{k} (-1)^k e^{-k\mu(t-\nu)}\]
which leads to 
\begin{align}\label{eq:pij2}
&p(i,j)= \nonumber\\
&\binom{c-1}{c-j-1} \frac{(c\mu)^{i-c+2}}{(i-c+1)!} \sum_{k=0}^{c-j-1} \binom{c-j-1}{k} (-1)^k \int_{0}^{\infty}\!\! e^{-\mu (j+k)t} \int_{0}^{t} \nu^{i-c+1} e^{-\mu(c-j-k)\nu} \,d\nu\,dA(t)\;.
\end{align}
Let ${\cal K} = \frac{(\mu(c-j-k))^{i-c+2}}{(i-c+1)!}$.
Now, noting that $ {\cal K}\nu^{i-c+1} e^{-\nu\mu(c-j-k)}d\nu $  is a gamma  density function function, we obtain ( e.g., Tijms\cite{Tij03}, page 442)
\begin{align*}
\int_{0}^{t}\nu^{i-c+1} e^{-\nu\mu(c-j-k)}d\nu
&= \frac{(i-c+1)!}{(\mu(c-j-k))^{i-c+2}} \bigg[1 - \sum_{r=0}^{i-c+1} \frac{(\mu(c-j-k))^r t^r e^{-t\mu(c-j-k)}} {r!} \bigg].\\
\end{align*}
Substituting  into (\ref{eq:pij2}), we have
\begin{align*}
p(i,j) &= 
\binom{c-1}{c-j-1} \frac{(c\mu)^{i-c+2}}{(i-c+1)!} \sum_{k=0}^{c-j-1} \binom{c-j-1}{k} (-1)^k \int_{0}^{\infty}\!\! e^{-\mu (j+k)t} \frac{(i-c+1)!}{(\mu(c-j-k))^{i-c+2}} \\
&\times \bigg[1 - \sum_{r=0}^{i-c+1} \frac{(\mu(c-j-k))^r t^r e^{-t\mu(c-j-k)}} {r!} \bigg]\,dA(t)\\
&= 
\binom{c-1}{c-j-1} \frac{(c\mu)^{i-c+2}}{(i-c+1)!} \sum_{k=0}^{c-j-1} \binom{c-j-1}{k} (-1)^k \frac{(i-c+1)!}{(\mu(c-j-k))^{i-c+2}} \\&\times \bigg[ \int_{0}^{\infty}\!\!  
e^{-\mu (j+k)t} \,dA(t) - \sum_{r=0}^{i-c+1} \frac{(\mu(c-j-k))^r}{r!} \int_{0}^{\infty}  t^r e^{-c\mu t} \,dA(t) \bigg]\;.
\end{align*}
Using, 
\[A^*(s) = \int_{0}^{\infty} e^{-st} dA(t)  \;\;\;\hbox{and  }\;
A^*_n(s) = \int_{0}^{\infty} t^n e^{-st}dA(t);\]
 we have
\begin{align*}
p(i,j) &= 
\binom{c-1}{c-j-1} \frac{(c\mu)^{i-c+2}}{(i-c+1)!} \sum_{k=0}^{c-j-1} \binom{c-j-1}{k} (-1)^k \frac{(i-c+1)!}{(\mu(c-j-k))^{i-c+2}} \\&\times \bigg[A^*(\mu (j+k)) - \sum_{r=0}^{i-c+1} \frac{(\mu(c-j-k))^r}{r!} A^*_r(c\mu) \bigg]\\
&= 
\binom{c-1}{c-j-1} \sum_{k=0}^{c-j-1} \binom{c-j-1}{k} (-1)^k \left(\frac{c}{c-j-k}\right)^{i-c+2} \\&\times \bigg[A^*(\mu (j+k)) - \sum_{r=0}^{i-c+1} \frac{(\mu(c-j-k))^r}{r!} A^*_r(c\mu) \bigg]\;.
\end{align*}
By a change of variable where  $c-j-k$ is replaced by $k$,  we obtain
\begin{align*}
p(i,j) &= 
\binom{c-1}{c-j-1} \sum_{k=1}^{c-j} \binom{c-j-1}{c-j-k} (-1)^{c-j-k} \left(\frac{c}{k}\right)^{i-c+2} \\&\times \bigg[A^*(\mu (c-k)) - \sum_{r=0}^{i-c+1} \frac{(k\mu)^r}{r!} A^*_r(c\mu) \bigg]\;.
\end{align*}
Expanding the combinations and simplifying leads to
\begin{align*}
p(i,j) &= 
\frac{(c-1)!}{(c-j-1)!(c-1-(c-j-1))!} \sum_{k=1}^{c-j} \frac{(c-j-1)!}{(c-j-k)!(c-j-1-(c-j-k))!} 
\\&\times (-1)^{c-j-k} \left(\frac{c}{k}\right)^{i-c+2} \bigg[A^*(\mu (c-k)) - \sum_{r=0}^{i-c+1} \frac{(k\mu)^r}{r!} A^*_r(c\mu) \bigg]\;.\\
\end{align*}
Rearrange to obtain
\begin{align}% \label{eq:direct}
p(i,j) =\sum_{k=1}^{c-j} \frac{(-1)^{c-j-k}(c-1)!}{(k-1)!(c-j-k)!j!}  \left(\frac{c}{k}\right)^{i-c+2} \bigg[A^*((c-k)\mu ) - \sum_{r=0}^{i-c+1} \frac{(k\mu)^r A^*_r(c\mu)}{r!}  \bigg].
\end{align}
which is the desired result.
\done
 
Using Lemma~\ref{lem:fund}, we obtain the main computationally useful result.

\begin{thm}\label{thm:pij}
(i) For $ j\leq i+1\leq c$, the transition probabilities of $p(i,j)$ are given by 
\[
p(i,j)=  \binom{i+1}{i-j+1}\sum_{r=0}^{i-j+1} (-1)^r \binom{i-j+1}{r} A^*((j+r)\mu))  \;. 
\]
(ii) For $c\leq j\leq i+1$, the transition probabilities of $p(i,j)$ are given by 
\[
p(i,j)=  \frac {(c\mu )^{i-j+1}A^*_{i-j+1}(c\mu)}{(i-j+1)!}\;.
\]
(iii) for  $j >0, j+1 \leq c\leq i$
\begin{eqnarray}
p(i,j)&=& \sum^{c-j}_{k=1} \frac{(-1)^{c-j-k}(c-k)C^a_{k, c-j}}{j}\left(\frac{c}{k}\right)^{i-c+2} 
 \left[ A^*((c-k)\mu)- \sum^{i-c+1}_{r=0} \frac{(k\mu )^{r}A^*_r(c\mu)}{r!}  \right] \;; 
\end{eqnarray}
where $C^a_{k, c-j}=\prod_{m=1}^{k-1} \frac{c-m}{k-m}\times \prod_{m=k+1}^{c-j} \frac{c-m}{m-k}\;,$  $\prod$  over empty sets is $1$ ;and $p(i,j)= 0$, otherwise.
\end{thm}
{\bf Proof.} One can easily see that (i) and (ii) of the Theorem  follow  from Lemma~\ref{lem:pij} (i) and (ii). 
%and equation~(\ref{eq:aij}) 
%one can easily  prove (i) and (ii) of the Theorem respectively. 
Part (iii) follows  from Lemma~\ref{lem:fund} by noting that
  after some simplification 
\begin{align}\label{exp:1}
C^a_{k, c-j}= \frac{(c-1)!(c-k-1)!}{(c-k)!(k-1)!(j-1)!(c-j-k)!} \;.
\end{align}
This completes the proof. \done

Similar results are given by El-Taha~\cite{Elt21} using a convolution approach. Our approach here is more direct and intuitive. We note that Lemma~\ref{lem:fund} allows for  direct evaluation of $p(i,j)$ for $j=0$ in region 3, while Theorem~\ref{thm:pij}(iii) does not. However we can always compute $p(i,0)$ as the complement of the remaining  $p(i,j)$ values for $j \geq 1$. 

\begin{figure}[ht] 
 \includegraphics[bb=15 10 600 400,width=6.00in,height=4.00in,keepaspectratio]{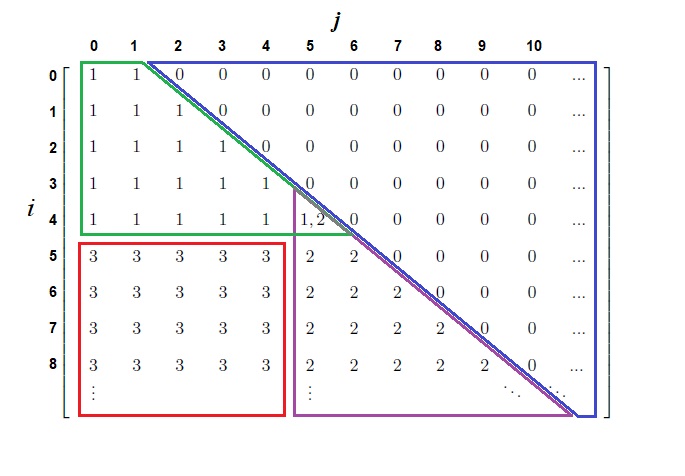}
  \caption{Transition Matrix Regions for  $c=5$}
 \label{FigA3}
\end{figure}

By referring to Figure~\ref{FigA3}, region $1$ is  where $p(i,j)$  is given by Lemma~\ref{lem:pij} (i) and Theorem~\ref{thm:pij} which is given by a binomial distribution function, region $2$ is where $p(i,j)$  is similar to the $GI/M/1 $ model with service rate $c\mu$,  region $3$ is  where $p(i,j)$  is given by
 the direct method, and region $4$ is where $p(i,j)$ is $0$.
\bigskip
   
Computing the one-step transition probabilities $p(i,j)$ requires the evaluation of the derivatives  of the $LST$ of the inter-arrival time distribution functions. This can be easily done exactly for  all Phase ($PH$) Type distribution functions. Because $PH$ Type distributions are dense in the set of all continuous distribution functions with support on $[0,\infty)$, our approach works at least approximately for any continuous non-negative distribution function. Now, we introduce $PH$ Type distribution functions and give their properties, see Neuts~\cite{Neu81} Ch 2, including a formula for the derivatives of their $LST$.

Phase-Type distribution functions are modeled as the time until absorption in a Markov process (MP) with a single absorption state. Consider a $PH$ Type distribution function  defined on a $MP$ with states $\{1,\cdots,k+1\}$, $k+1$ being the absorption state, and a transition rate matrix
\begin{equation*}
Q = 
\begin{pmatrix}
T& T^0 \\
\bf{0} & 0 
\end{pmatrix}, \mbox{where}
\end{equation*}
$T$ is a $k\times k$ matrix with $T_{ij} \geq 0, i\not= j,$ and $T_{ii} <0, i=1, \ldots, k$; $\bf{1}$ is a vector of ones; $T^0$ is chosen such that $T \bf{1}$$+T^0=\bf{0}$ i.e., the row of elements of $Q$ add up to $0$; 
and ($\bos{\alpha}, 0$) with ${\bos{\alpha} \bf{1}}=1$,
 is the initial distribution of the MP. Let the random variable $X$ be the time until absorption, then the distribution function of $X$ is said to be a $PH$-Type distribution function with representation (${\bos \alpha}, T$). Now let the distribution function of $X$ be the inter-arrival time distribution function of the  $GI/M/c$ model, then the distribution function and its $LST$ are given by

\begin{eqnarray*}
A(t)&=&1-\bos{\alpha}\exp({Tt}){\bf 1}, t\geq 0;\\
A^*(s) &=& \bos{\alpha}(sI-T)^{-1}T^0 \;.
\end{eqnarray*}
Moreover, 

\begin{equation}\label{eq:PHLST}
A^*_n(s) = n!\bos{\alpha}(sI-T)^{-(n+1)}T^0 \;.
\end{equation}

This is explored further in Section~\ref{sec:app} where we give closed form expressions for the derivatives of several cases of 
$PH$ Type distribution functions including the exponential, the Erlang, and the hyper-exponential.
\bigskip

\noindent
{\bf Remarks on Complexity and Stability.}
We note that the method given by Theorem~\ref{thm:pij} (iii) involves two finite sums and computations of the order of  $O((i-c)(c-j))$.  In contrast the method using Lemma~\ref{lem:pij} (iii) would involve numerically approximating  a double infinite integration. El-Taha~\cite{Elt21} also discusses  numerical stability issues that arise  when calculations involve  alternating between positive and negative terms. Specifically, he  suggests collecting all positive (negative terms) together and using one subtraction at the end.

We point out that the contribution of Lemma~\ref{lem:fund} and Theorem~\ref{thm:pij} is in computing the transition probabilities efficiently using a finite sum as compared to numerical integration that would be needed in part (iii) of Theorem~\ref{thm:pij}. In the next two subsections we provide efficient methods to compute the arrival-time and time-average probabilities. 

\subsection{Arrival-Times Stationary Distribution}\label{subsec:ATSD}
%\bigskip

Let $\{\pi(i), i= 0, \ldots, \}$ be the stationary distribution function of the imbedded Markov chain
$\{X_n, n=0,1\ldots\}$, so that 
$
\pi(k) = \lim_{n \rightarrow~\infty} P \{ X_n=k\}\;. 
$ 
That is, $\pi(k)$ is the probability an arrival sees $k$ customers.
We assume a stable system which is the case when $\rho=\lambda/c\mu <1$ in the infinite capacity model. We use the concept that probability flow across cuts balance to efficiently solve the pre-arrival time probabilities for the $GI/M/c/N$ model. For $j=1,2,\cdots,N$ (using Kelly~\cite{Kel79}, Lemma 1.4), we obtain
%
%
%
%
%We give  a recursive  algorithm to solve for the departure-time probabilities ${\pi(.)}$ for the $GI/M/c/N$ model.
%Using the  fact that probability flow across cuts balances (Kelly~\cite{Kel79}, Lemma 1.4) one can  show that for $j=1,2,\cdots,N$ 
%
\begin{eqnarray*}
\pi(j)p(j,j+1)&= & \sum_{k=j+1}^{N} \pi(k)\sum_{i=0}^{j}p(k,i)
\end{eqnarray*}
which can be rewritten as 
\begin{equation}\label{eq:pi_j+1}
\pi(j)=  \sum_{k=j+1}^{N} \pi(k)a(k,j)/p(j,j+1)
\end{equation}
where $a(k,j)=\sum_{i=0}^{j}p(k,i)$.
We note that solving (\ref{eq:pi_j+1}) is more efficient and numerically stable than solving the corresponding global balance equations. This was demonstrated  by Stidham~\cite{Sti87}.

Using (\ref{eq:pi_j+1}) we compute $\pi(.)$ interativelly starting with $\pi(N)$. For the infinite capacity model we can iterate on $N$ until we achieve a desired level of accuracy. Instead, we pre-determine $N$ so that  a prescribed level of precision is achieved.
It is well-known (Gross and Harris~\cite{Gro08}) that $\pi(n) =C \sigma^{n}$ for $n \geq c$  where $C$ is a constant 
and
$|\sigma| <1$ is the solution of the  equation
 \begin{equation} \label{eq:sigma}
\sigma =\int_0^{\infty} e^{-c\mu(1-\sigma)}t dA(t)\equiv A^*(c\mu(1-\sigma))\;.
\end{equation} 
Moreover,  $C$ is  given by Gross and Harris~\cite{Gro08} as
\[
C= \frac{1-\sum_{j=0}^{c-1}\pi(j)}{\sigma^c (1-\sigma)^{-1}}\;.
\]
The following lemma serves to guide our selection of $N$.
\begin{lem}\label{lem:n}
If we choose 
\[
n > c+  \frac{\ln (\epsilon) - 2 \ln (1-\sigma)}{\ln (\sigma)}\;;
\]
then $n$ will  satisfy  the  criterion $|\pi(n+1)-\pi(n)|< \epsilon$, where  $\epsilon$ is the tolerance level.
\end{lem}
\noindent{\bf Proof.}  It follows from $|\pi(n+1)-\pi(n)|< \epsilon$ that 
\begin{eqnarray*}
\sigma^n&< & \epsilon/(K(1-\sigma))\\
 &=& \epsilon \sigma^c/((1-\sigma)^2 (1-\sum_{j=0}^{c-1}\pi(j)))
\end{eqnarray*}
 Take the natural logarithm of both sides to obtain
\begin{eqnarray*}
n&> & \frac{\ln(\epsilon) + c \ln(\sigma) - 2\ln(1-\sigma) - \ln(1-\sum_{j=0}^{c-1}\pi(j))}{\ln(\sigma)}\\
 &>&  c+ \frac{\ln(\epsilon) - 2\ln(1-\sigma)}{\ln(\sigma)}  - \frac{\ln(1-\sum_{j=0}^{c-1}\pi(j))}{\ln(\sigma)}\;.
\end{eqnarray*}
This completes the proof.\done

Computing the  stationary probabilities $\{\pi(j); j \geq 0\}$ involves the following steps:
 initialize the system by determining precision $\epsilon$, number of servers $c$,  the inter-arrival times distribution 
function, the buffer $N$ (Lemma~\ref{lem:n}), and computing $\sigma$ (using \ref{eq:sigma})\;; compute transition probabilities $p(i,j)$ (Theorem~\ref{thm:pij}) and $a(k,j)$;  compute  $\pi'(j)=\sigma^j$\; ($j\geq c$);  compute  $\pi'(j)=\sum_{k=j+1}^N\pi'(k)a(k,j)/p(j,j+1)$\;($j\leq c$); and finally normalize.
 In the Appendix, we give a detailed algorithm to compute  the imbedded and time average probabilities.

\subsection{Time-Average Stationary Distribution}

One is typically interested in the system performance measures like the mean number of customers in the system and the queue $L$ and $L_q$, and  the mean delay in the system and the queue $W$ and $W_q$ respectively. To do this we need the time-average probability distribution function which is of interest in its own right. Let $\{X(t), t\geq \}$ be a stochastic process so that $X(t)$ is the number of customers in the system at time $t$, then the  
 time-average stationary distribution is defined as 
 \[
p(n) = \lim_{t \rightarrow~\infty} P \{ X(t)=n \}, n=0,1, \ldots, .
\]
We evaluate the stationary distribution function $p(.)$ by relating those probabilities to the  pre-arrival probabilities $\pi(.)$. We know that (El-Taha and Stidham~\cite{Elt99})
\[ 
\lambda \pi(n)=\lambda_n p(n)=\mu_{n+1}p(n+1)\;.
\]
Using this relation leads to the following result.
\begin{lem}\label{lem:p(i)}  Consider the $GI/M/c/N$ queueing model and let $\rho=\lambda/c\mu$.
Then the  system size probabilities, $p(i), i=0, \ldots N$, are given by

\[
p(n)= \left\{ \begin{array}{lll} 
                (1-\rho) + \rho \pi(N)- \rho\sum_{k=0}^{c-2}\frac{c-k-1}{k+1}\pi(k), &  n=0\;;\\

                c\rho \pi(n-1)/n\;,            & 1\leq n\leq c \;;\\
                 \rho \pi(n-1) \;,         & c< n \leq N \;.
                \end{array}
                \right. \\
\]  
\end{lem}
A full proof of Lemma~\ref{lem:p(i)} is given by  El-Taha~\cite{Elt21}.
For the stable infinite buffer $GI/M/c/\infty$ model, with $\rho<1$,  $p(0)$ is  computed as
\begin{equation}
p(0)= (1-\rho)- \rho\sum_{k=0}^{c-2}\frac{c-k-1}{k+1}\pi(k)\;.
\end{equation}

%{\bf Proof.}
%We know that (El-Taha and Stidham~\cite{Elt99})
%\[ 
%\lambda \pi(n)=\lambda_n p(n)=\mu_{n+1}p(n+1)\;.
%\]
%Therefore,
%\[ 
% p(n)= \lambda \pi(n-1)/\mu_{n}\;; n\geq 1\;.
%\]
%Thus, $p(n)= c\rho \pi(n-1)/n$  for $1\leq n\leq c$, and  $p(n)= \rho \pi(n-1)$ for  $n> c$. We still need to compute $p(0)$.
%Normalize  the $p(i)'s$ to obtain
%\begin{eqnarray*}
%1-p(0) &=& \sum_{k=1}^N p(k)\\
%&=& \sum_{k=1}^{c} c\rho \pi(k-1)/k  +\sum_{k=c+1}^{N} \rho \pi(k-1)\\
%   &=& c\rho\sum_{k=0}^{c-1}  \pi(k)/(k+1)  + \rho [1- \sum_{k=0}^{c-1}  \pi(k)-\pi(N)]\;.
%\end{eqnarray*}
%Therefore
%%
%\begin{eqnarray*}
%p(0) &=&  (1-\rho) - c\rho\sum_{k=0}^{c-1}  \pi(k)/(k+1) +  \rho  \sum_{k=0}^{c-1}  \pi(k) +\rho \pi(N)\\
%&=&  (1-\rho) - \sum_{k=0}^{c-1}(\frac{c\rho}{k+1}-\rho)  \pi(k)  + \rho \pi(N)\\
%   &=&   (1-\rho) - \rho \sum_{k=0}^{c-2}\frac{c-k-1}{k+1}  \pi(k)+\rho \pi(N)\;;\\
%\end{eqnarray*}
%which completes the proof.\done
%
%\noindent
%{\bf Remark.} For the stable infinite buffer model $GI/M/c/\infty$ model, with $\rho<1$,  $p(0)$ is  computed as
%
%\begin{equation}
%p(0)= (1-\rho)- \rho\sum_{k=0}^{c-2}\frac{c-k-1}{k+1}\pi(k)\;.
%\end{equation}
%
%
%With  $p(i)$'s determined from Lemma~\ref{lem:p(i)}, other measures of performance are now determined
%easily. We focus on a few measures of performance that we believe best reflect the overall performance of the system. 
%The mean number of  customers  $L$ is given by
%$
%E[L]=\sum_{i=1}^{\infty}ip(i)\;.
%$
%Little's law implies  that the mean delay in the system
% is  
%$
%W=L/\lambda\;.
%$
In  Section~\ref{sec:app} we calculate $L$, $W$, $L_q$ (the mean number of customers in the queue), and  $W_q$ (the mean delay in the queue)
using several  service time distributions.

\section{Extension to $GI^X/M/c/N$ Model}\label{sec:extensions}

In this section we discuss an intuitive method to extend our $GI/M/c$ results to the $GI^X/M/c$ and  $GI^X/M/c/N$ models.  We note that the models have been discussed by several articles, notably by  Laxmi and Gupta~\cite{Lax00} and Chaudry and Kim~\cite{Cha16} among others mentioned in the introduction. Our approach is different in the sense that we generate the one-step transition probabilities of the batch arrival models directly from the standard $GI/M/c$ model and the batch size distribution function.

The section builds on Theorem~\ref{thm:pij} by devising a remarkably intuitive and new approach to express the transition probabilities for the batch arrival $GI/M/c$ models in terms of the transition probabilities of the of the $GI/M/c$ model given by Theorem~\ref{thm:pij}. These results are expressed in Lemmas~\ref{lem:pij-batch}, \ref{lem:pij-partial}, and \ref{lem:pij-full}. Also the results in Lemmas~\ref{lem:p(i)-batch}, \ref{lem:p(i)-partial}, and \ref{lem:p(i)-full}, that follow from more general results in El-Taha and Stidham~\cite{Elt99} allow us to compute time-average  from arrival-time probabilities. 

Let $X$ be a random variable that represents the batch size, and let $b(k)$ be the $pmf$ of $X$  with support $[1, \infty)$ such that  $b(k)=P(X=k)$. We note here that  typical random variables that are used  to model batch sizes have a support that starts at $0$. Let $g(k)$ be a $p.m.f.$ with support $[0,\infty)$, that is, with g(k) we allow a batch of size $0$ to occur like the Poisson and the geometric random variables. In this situation, we think of $b(k)$ as $b(k)=g(k)/(1-g(0))$.  Also note that in these models $A(t)$ is the distribution function of batch inter-arrival times and $1/\lambda$ is the mean time between batch arrivals, so that the overall arrival rate $\lambda_A = \lambda E[X]$. Also  $\rho=\lambda E[X]/c\mu$. We require $\rho <1$ for the infinite buffer models.

Moreover let $p^*(i,j)$ be the one step transition probabilities  just before arrival of the imbedded Markov chain of the $GI^X/M/c$ and $GI^X/M/c/N$ models. Recall that $p(i,j)$ are the corresponding one step transition probabilities associated with the $GI/M/c$ and $GI/M/c/N$  models. 
Now we consider three cases.
\medskip

\noindent
{\bf Case 1.  $GI^X/M/c$}

We consider the  $GI/M/c$ with batch arrivals and obtain the one step-transition probabilities.

\begin{lem}\label{lem:pij-batch}
Consider the $GI^X/M/c$ queueing model. Then, for all $i,j$
\[  p^*(i,j)= \sum_{k=1}^{\infty} p(i+k-1,j)b(k)\;.
\]
\end{lem}

Now replace $p(i,j)$ with $p^*(i,j)$ to compute $\{\pi(i)\}$ of the resulting Markov chain.
In order to relate the time-average to pre-arrival probabilities probabilities we use the following  result. 

\begin{lem}\label{lem:p(i)-batch}  Consider the $GI^X/M/c$ queueing model with batch arrivals. 
Then the  system size probabilities, $\{p(i), i=0, \ldots \}$, are given by
\[
p(n)= \left\{ \begin{array}{lll} 
      \lambda\sum_{k=1}^{n-1}\pi(k)B^c(n-k)/\min(n,c)\mu & n=1,\cdots \;;\\ \\
       1-\sum_{k=1}^{\infty}\pi(k) \;,         & n=0 \;;
                \end{array}
                \right. \\
\] 
 where $B^c(n-k)=\sum_{n-k}^{\infty}b(i)$.
\end{lem}
{\bf Proof.} The proof follows from  El-Taha and Stidham~\cite{Elt99}, equation (4.34), page 107 ; see also Stidham and El-Taha~\cite{Sti89}. 
\done

\medskip
\noindent
{\bf Case 2: $GI^X/M/c/N$ with partial rejection.}

In this case a batch that brings the system state above $N$ is partially accepted  in the sense that the system will accept a part of the batch for the state to reach $N$ and the rest of the batch is rejected. Then  for all $i,j$ the one step-transition probabilities  are given by this result.

\begin{lem}\label{lem:pij-partial}
Consider the $GI^X/M/c/N$ queueing model with partial rejection. Then, for all $i,j$
\[  p^*(i,j)= \sum_{k=1}^{N-i-1} p(i+k-1,j)b(k) + p(N-1,j)B^c(N-i)\;;
\]
where $B^c(N-i)=\sum_{k=N-i}^{\infty}b(k)$.
\end{lem}
{\bf Proof.} It is clear that 
\begin{eqnarray*}
p^*(i,j)&= &\sum_{k=1}^{N-i-1} p(i+k-1,j)b(k) + p(N-1,j)(\sum_{k=N-i}^{\infty}b(k))\\
    &= &\sum_{k=1}^{N-i-1} p(i+k-1,j)b(k) + p(N-1,j)B^c(N-i)\;.
\end{eqnarray*}

In order to relate the time-average to pre-arrival probabilities probabilities we use the following  result. 

\begin{lem}\label{lem:p(i)-partial}  Consider the $GI^X/M/c/N$ queueing model with partial rejection. 
Then the  system size probabilities, $p(i), i=0, \ldots N$, are given by
\[
p(n)= \left\{ \begin{array}{lll} 
      \lambda\sum_{k=1}^{n-1}\pi(k)B^c(n-k)/\min(n,c)\mu, & n=1,\cdots N\;;\\ \\
       1-\sum_{k=1}^N\pi(k) \;,         & n=0 \;;
                \end{array}
                \right. \\
\] 
 where $B^c(n-k)=\sum_{n-k}^{\infty}b(i)$.
\end{lem}
{\bf Proof.} The proof is similar to Lemma~\ref{lem:p(i)-batch}. 
\done

\medskip 

\noindent
{\bf Case 3: $GI^X/M/c/N$ with full rejection.}

In this case a batch that brings the system state above $N$ is fully rejected. 

\begin{lem}\label{lem:pij-full}
 Consider the $GI^X/M/c/N$ with full rejection. The, for all $i,j$
\[  p^*(i,j)= \sum_{k=1}^{N-i} p(i+k-1,j)b(k) + p(i,j)B^c(N-i+1)\;;
\]
where $B^c(N-i+1)=\sum_{k=N-i+1}^{\infty}b(k)$.
\end{lem}
{\bf Proof.}
It is clear that 
\begin{eqnarray*}
p^*(i,j)&= &\sum_{k=1}^{N-i} p(i+k-1,j)b(k) + p(i,j)(\sum_{k=N-i+1}^{\infty}b(k))\\
    &= &\sum_{k=1}^{N-i-1} p(i+k-1,j)b(k) + p(i,j)B^c(N-i+1)\;.
\end{eqnarray*}

 In order to relate the time-average probabilities to arrival epochs  probabilities we use the following  result. 

\begin{lem}\label{lem:p(i)-full}  Consider the $GI^X/M/c/N$ with full rejection.
Then the  system size probabilities, $p(i), i=0, \ldots N$, are given by
\[
p(n)= \left\{ \begin{array}{lll} 
      \lambda\sum_{k=1}^{n-1}\pi(k)[B^c(n-k)-B^c(N-k+)]/\min(n,c)\mu & n=1,\cdots N\;;\\\\
       1-\sum_{k=1}^N\pi(k) \;,         & n=0 \;;
                \end{array}
                \right. \\
\] 
 where $B^c(n-k)=\sum_{n-k}^{\infty}b(i)$.
\end{lem}
{\bf Proof.} The proof is similar to Lemma~\ref{lem:p(i)-partial}, except here  batches that result in more than $N$  customers in the system are  totally rejected. \done 

Yao et.al.\cite{Yao84} and Laximi and Gupta~\cite{Lax00} give similar relations to relate  time-average and pre-arrival probabilities as in Lemma ~\ref{lem:p(i)-partial} and Lemma~\ref{lem:p(i)-full}.
  
 Replacing $p(i,j)$ with $p^*(i,j)$, the arrival-time probabilities can be computed using $\pi=\pi P^*$, $\sum_{i\in S} \pi(i)=1$ where $P^*=[p^*(ij)]$ is the one-step transition matrix and $S$ is the state space. One can also use  one of the methods discussed in the literature in the introduction to compute the arrival-time distribution.
The time-average probabilities are then computed using the corresponding Lemma ~\ref{lem:p(i)-batch}, Lemma ~\ref{lem:p(i)-partial}, or Lemma~\ref{lem:p(i)-full}. In Appendix A, the computations for the batch arrival models are given by an appropriate adaptation of the algorithm  for the $GI/M/c/N$ model.

%%%%%%%%%%%%%%%%%%%%%%%%%%%%%%%%%%%%%%%%%%%%%%%%%%%%%%%%%%%%%%%%%%%%%%%%%%%%%%%%%%%%%%%%%%%%%%%%%%%%%
\section{Applications } \label{sec:app}
%%%%%%%%%%%%%%%%%%%%%%%%%%%%%%%%%%%%%%%%%%%%%%%%%%%%%%%%%%%%%%%%%%%%%%%%%%%%%%%%%%%%%%%%%%%%%%%%%%%%%

In this section we consider both finite buffer (small scale) and  infinite buffer (large scale) examples of the $GI/M/c$ model. In the small scale case we focus on small finite buffers versus infinite buffers for the large scale case. We note that in the infinite  buffer case we truncate the system size  by using $|\pi(n+1)-\pi(n)|< \epsilon$ and appealing to Lemma~\ref{lem:n}.
We start with small scale finite buffer examples in  the next subsection. 

For both small scale and large scale applications  we focus on inter-arrival time distributions whose  $LST$ have closed form multiple derivatives. We choose four distribution functions with coefficients of variation that vary from $0$ to infinity. Specifically, we  select the deterministic, Erlang, exponential, and the hyper-exponential distribution functions. Let $a(t)$ be the $p.d.f.$ of the inter-arrival times.
Recall that we require that the  inter-arrival time distribution function has a mean $E[A]=1/\lambda$. Here we give explicit forms for the derivatives of the $LST$ of the inter-arrival time distribution functions.
\medskip

\noindent
{\bf Deterministic.} In this case we assume that  $a(t)=a$ w.p. 1, so that
$
A^*_n(s)= a^ne^{-sa}\;.
$
 Note that  here $\lambda =1/a$.
\medskip

\noindent
{\bf Exponential.} Here $a(t)= \lambda e^{-\lambda t}, t \geq 0$, so that
$
A^*_n(s)= n!\lambda/(s+ \lambda)^{n+1}\;.
$
\medskip

\noindent
{\bf Erlang.} The density function for a $k$ phase Erlang is  $a(t)= \frac{\theta(\theta t)^{k-1}}{(k-1)!}e^{-\theta t}, \theta >0,t \geq 0$. Here
\medskip
$
A^*_n(s)= n!\binom{k+n-1}{k-1}\theta^k/(s+ \theta)^{k+n}\;.
$
Note that here $E[A]=k/\theta$, so that $\theta=k\lambda$.
We would like the mean to stay constant. So, we replace $\theta$ by $k\lambda$. The density function for a $k$ phase Erlang is  $a(t)= \frac{k\lambda(k\lambda t)^{k-1}}{(k-1)!}e^{-k\lambda t}, \lambda >0,t \geq 0$. Here
\medskip
$
A^*_n(s)= n!\binom{k+n-1}{k-1}(k\lambda)^k/(s+ k\lambda)^{k+n}\;.
$
Note that again $E[A]=1/\lambda$.

\noindent
{\bf Hyper-exponential.} Let $a_i(t)$ be an exponential $p.d.f$ with parameter $\lambda_i$. Then, the $k$ phase hyper-exponential is given by 
$
 a(t) = \sum_{i=1}^{k} p_i a_i(t);
$
 so that $ A^*(s) = \sum_{i=1}^{k} \frac{p_i\lambda_i}{s+\lambda_i}$.
We shall use the two phase hyper-exponential  which is   a mixture of two
exponential distribution functions. The density function can be written as  
$a(t)= p\lambda_1 e^{-\lambda_1 t}+(1-p)\lambda_2 e^{-\lambda_2 t}, \lambda_1 >0, \lambda_2 >0, t \geq 0, (0 \leq p \leq 1$). 
Therefore
\[
A^*_n(s)= n! \left[ \frac{p\lambda_1}{(s+\lambda_1)^{n+1}} + \frac{(1-p)\lambda_2}{(s+\lambda_2)^{n+1}}\right]\;.
\]
In the next subsection we show how to use these transforms to write  the one step transition matrix in explicit form.

\subsection{Small Scale Applications}

Here, we present small-scale, finite-buffer examples for a system with $c=3$ servers and a total capacity of $K=6$.
Now referring to  the regions  described in Figure \ref{FigA3} and using Theorem~\ref{thm:pij}, we write  the one-step transition probabilities in the more computationally convenient form as follows.
\begin{enumerate}[label=(\roman*)]   
	\item For Region 1, where $i\leq c-1$ and $j\leq i+1$, (i.e.: $i=0,1,2 $ and $j=1,...,i+1$), we use
	\begin{align*}p(i,j) &= \frac{(i+1)!}{j!} \sum_{k=0}^{i-j+1} \frac{(-1)^k A^*((j+k)\mu)}{(i-j-k+1)!k!}\;.
	\end{align*}
	
	\item For Region 2, where $i=c,c+1,...,K-1$,\;\, $j=c,c+1,...,i+1, \;\;i+1\leq K \;\;$ (i.e.: $i=3,4,5$ and $j=3,4,5$), we use
\begin{align*}
p(i,j) &= \frac{(c\mu)^{i-j+1}}{(i-j+1)!} A^*_{i-j+1}(c\mu) 
= \frac {(3\mu )^{i-j+1}A^*_{i-j+1}(3\mu)}{(i-j+1)!}\;.
\end{align*}

	\item For Region 3, where $1\leq j\leq c-1 < i$ (i.e.: $i = 3, 4, 5$ and $j = 1, 2$),
	we use the direct result 
	\[p(i,j) = \sum_{k=1}^{3-j} \frac{C_{k,3-j}(3-k)}{j}\bigg(\frac{3}{k}\bigg)^{i-1} \Bigg[A^*((3-k)\mu) - \sum_{r=0}^{i-2}\frac{(k\mu)^r A^*_r(3\mu)}{r!} \Bigg] \;; \]
	\par
	where
 \[C_{k,3-j} = \prod_{m=1}^{k-1}\frac{3-m}{k-m}\prod_{m=k+1}^{3-j}\frac{3-m}{k-m} \;\;\; \text{, } \; C_{3-j,3-j} = \prod_{m=1}^{2-j} \frac{3-m}{3-j-m} \;,\;\; \text{ and }\] \[ \; C_{v,v}=\prod_{m=1}^{v-1}\frac{3-m}{v-m} \;;
	\text{ so that } C_{1,1}=1, C_{1,2}=-1, C_{2,2}=2\;.\]
\end{enumerate}

	Also,  for $i\geq c$ and $j=0$, we simply use
	$p(i,j) = 1 - \sum_{n=1}^{K}p(i,n)$. Moreover, for Region 4, $p(i,j) = 0$ where $j>i+1$. When $i=K$, use $p(K,j)=p(K-1,j)$ for all $j=0,\cdots,K$. Note that, if $i$ is the number in the system immediately prior to an arrival, then the  transition probabilities when $i=K-1$ must be the same as the probabilities when $i=K$, because in the first case, the system becomes full, and in the second case, the system is already full and the new arrival is lost.  Given the Markovian property, the  transition probabilities  are unaffected by additional arrivals while the system is full. 
Therefore the one-step transition matrix is given by
\bigskip

\tiny
$P = {
		\begin{bmatrix*}[l]
			1 - A^*(\mu) &\;\; A^*(\mu) &\;\; 0 &\;\; 0 &\;\; 0 &\;\; 0&\;\; 0\\[0.125in]
			
			1 - 2A^*(\mu) + A^*(2\mu) &\;\; 2A^*(\mu) - 2 A^*(2\mu) &\;\; A^*(2\mu) &\;\; 0 &\;\; 0 &\;\; 0 &\;\; 0\\[0.125in]
			
			1 - 3A^*(\mu) + 3A^*(2\mu) &\;\; 3A^*(\mu) - 6A^*(2\mu) &\;\; 3A^*(2\mu) - 3A^*(3\mu) &\;\; A^*(3\mu) &\;\; 0 &\;\; 0 &\;\; 0\\[0.01in]
			  
                        \;\;\;\;- A^*(3\mu) &\;\;\;\;\; + 3A^*(3\mu) &&&&& \\[0.125in]
			1 - \sum_{n=1}^{4}p(3,n) 
			& p(3,1) & p(3,2) &\;\; 3\mu A^*_1(3\mu) &\;\; A^*(3\mu) &\;\; 0 &\;\; 0  \\[0.125in]
			
			1 - \sum_{n=1}^{5}p(4,n) 
			& p(4,1) & p(4,2) &\;\; \frac{9}{2}\mu^2 A^*_2(3\mu) &\;\; 3\mu A^*_1(3\mu) &\;\; A^*(3\mu) &\;\; 0 
			\\[0.125in]
			
			1 - \sum_{n=1}^{6}p(5,n) 
			& p(5,1) & p(5,2) &\;\; \frac{9}{2}\mu^3 A^*_3(3\mu) &\;\; \frac{9}{2}\mu^2 A^*_2(3\mu) &\;\; 3\mu A^*_1(3\mu) &\;\; A^*(3\mu)
			\\[0.125in]
			
			1 - \sum_{n=1}^{6}p(6,n) 
			& p(6,1) & p(6,2) &\;\; \frac{9}{2}\mu^3 A^*_3(3\mu) &\;\; \frac{9}{2}\mu^2 A^*_2(3\mu) &\;\; 3\mu A^*_1(3\mu) &\;\; A^*(3\mu)\\[0.03in]
			
	\end{bmatrix*} } 
\\
$
\normalsize

\noindent where 
\ls{1.2}
\begin{align*}
p(3,1) &= \frac{9}{2}A^*(\mu)- 18A^*(2\mu) + \frac{27}{2}A^*(3\mu) + 9\mu A^*_1(3\mu)\;;\\
p(3,2) &= 9A^*(2\mu) - 9A^*(3\mu) -9\mu A^*_1(3\mu)\;; \\
p(4,1) &= \frac{27}{4}A^*(\mu) - 54A^*(2\mu) + \frac{189}{4}A^*(3\mu) + \frac{81}{2}\mu A^*_1(3\mu) + \frac{27}{2} \mu^2 A^*_2(3\mu)\;; \\
p(4,2) &= 27A^*(2\mu) - 27A^*(3\mu) -27\mu A^*_1(3\mu) - \frac{27}{2}\mu^2 A^*_2(3\mu)\;; \\
p(5,1) &= \frac{81}{8}A^*(\mu) - 162A^*(2\mu) + \frac{1215}{8}A^*(3\mu) + \frac{567}{4}\mu A^*_1(3\mu) + \frac{243}{4} \mu^2 A^*_2(3\mu)  + \frac{27}{2} \mu^3 A^*_3(3\mu)\;; \\
p(5,2) &= 81A^*(2\mu) - 81A^*(3\mu) - 81\mu A^*_1(3\mu) - \frac{81}{2}\mu^2 A^*_2(3\mu) - \frac{27}{2}\mu^3 A^*_3(3\mu)\;; \\
p(6,1) &= p(5,1)\;; \hbox{and  } 
p(6,2) = p(5,2)\;.
\end{align*}

\noindent
{\bf Deterministic Arrivals.} 
For deterministic  inter-arrivals with probability density function  $a(t)= a$, and $0$ otherwise (i.e. $\lambda^{-1} = a)$, thus we have
\[A^*_n(s) = a^n e^{-sa}\;; \]
therefore
\[A^*(\mu(j+m)) = e^{-\mu(j+m) a} \;;\]
and
\[A^*_n(c\mu) = a^n e^{-c\mu a} \;.\]
	This gives the following:
$ A^*(0) =1, 
	A^*(\mu)=e^{-a\mu},
	A^*_1(3\mu) =a e^{-3a\mu},
	A^*(2\mu) =e^{-2a\mu},
	A^*_2(3\mu)  =a^2 e^{-3a\mu},
	A^*(3\mu) =e^{-3a\mu},$ and $
	A^*_3(3\mu)  =a^3 e^{-3a\mu}\;.$ Substitute in the  general arrivals matrix to get the corresponding  one-step transition matrix for  deterministic inter-arrival times.
Thus our transition matrix is as follows:
\bigskip

\tiny
$P = {
\begin{bmatrix*}[l]
	1 - e^{-a\mu} &\;\; e^{-a\mu} &\;\; 0 &\;\; 0 &\;\; 0 &\;\; 0 &\;\; 0\\[0.125in]
	
	1 - 2e^{-a\mu} + e^{-2a\mu} &\;\; 2e^{-a\mu} - 2e^{-2a\mu} &\;\; e^{-2a\mu} &\;\; 0 &\;\; 0 &\;\; 0 &\;\; 0 \\[0.125in]
	
	1 - 3e^{-a\mu} + 3e^{-2a\mu}  &\;\; 3e^{-a\mu} - 6e^{-2a\mu}  &\;\; 3e^{-2a\mu} - 3e^{-3a\mu} &\;\; e^{-3a\mu} &\;\; 0 &\;\; 0 &\;\; 0 \\[0.01in]

	\;\;\;\;\;\; -\; e^{-3a\mu} &\;\;\;\;\ +\; 3e^{-3a\mu} &&&\;\;  &\;\;  &\;\;  \\[0.125in]
	
	1 - \sum_{n=1}^{4}p(3,n) 
	&\;\; p(3,1) &\;\; p(3,2) &\;\; 3a\mu e^{-3a\mu} &\;\; e^{-3a\mu} &\;\; 0  &\;\; 0 \\[0.125in]
	
	1 - \sum_{n=1}^{5}p(4,n) 
	&\;\; p(4,1) &\;\; p(4,2) &\;\; \frac{9}{2}a^2\mu^2 e^{-3a\mu} &\;\; 3a\mu e^{-3a\mu} &\;\; e^{-3a\mu} &\;\; 0 \\[0.125in]
	
	1 - \sum_{n=1}^{6}p(5,n) 
	&\;\; p(5,1) &\;\; p(5,2) &\;\; \frac{9}{2}a^3\mu^3 e^{-3a\mu} &\;\; \frac{9}{2}a^2\mu^2 e^{-3a\mu} &\;\; 3a\mu e^{-3a\mu} &\;\; e^{-3a\mu}
	\\[0.125in]

	1 - \sum_{n=1}^{6}p(6,n) 
	&\;\; p(6,1) &\;\; p(6,2) &\;\; \frac{9}{2}a^3\mu^3 e^{-3a\mu} &\;\; \frac{9}{2}a^2\mu^2 e^{-3a\mu} &\;\; 3a\mu e^{-3a\mu} &\;\; e^{-3a\mu}
	\\[0.02in]
	
	\end{bmatrix*} } $
\\
 \normalsize
\par

\noindent where 
\ls{1.2}
\begin{align*}
p(3,1) &= \frac{9}{2}e^{-a\mu} - 18e^{-2a\mu}+ \frac{27}{2}e^{-3a\mu} + 9\mu ae^{-3a\mu}\;; \\
p(3,2) &= 9e^{-2a\mu} - 9e^{-3a\mu} - 9\mu ae^{-3a\mu}\;; \\
p(4,1) &= \frac{27}{4}e^{-a\mu} -54e^{-2a\mu} + \frac{189}{4}e^{-3a\mu} + \frac{81}{2}\mu ae^{-3a\mu} + \frac{27}{2} \mu^2 a^2e^{-3a\mu} \;;\\
p(4,2) &= 27e^{-2a\mu} - 27e^{-3a\mu} - 27\mu ae^{-3a\mu} - \frac{27}{2}\mu^2 a^2e^{-3a\mu}\;; \\
p(5,1) &= \frac{81}{8}e^{-a\mu} - 162e^{-2a\mu} + \frac{1215}{8}e^{-3a\mu} + \frac{567}{4}\mu ae^{-3a\mu} + \frac{243}{4}\mu^2 a^2e^{-3a\mu}  + \frac{27}{2} \mu^3 a^3e^{-3a\mu}\;;\\
p(5,2) &= 81e^{-2a\mu} - 81e^{-3a\mu} -81\mu ae^{-3a\mu} - \frac{81}{2}\mu^2 a^2e^{-3a\mu} - \frac{27}{2}\mu^3 a^3e^{-3a\mu}\;;\\
p(6,1) &= p(5,1), \hbox{ and }\;\;
p(6,2) = p(5,2)\;.
\end{align*}
\ls{1.5}
Similar to the deterministic distribution function, one can  easily generate one-step transition matrices  for the other distributions, namely the exponential, the Erlang, and the hyper-exponential.
%%%%%%%%%%%%%
\subsection{Numerical Results}
%%%%%%%%%%%%%%%%%%

Here we provide numerical results for the distributions given above, for a system with $c=3$ servers, arrival rate $\lambda=5$, service rate $\mu=2$, capacity of $N=6$, and utilization factor $\rho=.83\overline{3}$. For the hyper-exponential  we use $p=.8, \lambda_1=8$, and $\lambda_2=2$, so that the overall $\lambda =5$. Numerical results are reported in Tables 1-3. In Table 1, we use the exponential distribution function and compare our method with  the traditional method, i.e. the $M/M/3/6$ model. The results  match perfectly. This helps  to verify our approach numerically. Table 2 reports the time-average probabilities for the other three distribution functions, and Table 3 reports performance measures.

%%%%
%%%% Table 1 
\begin{center}
\noindent
\begin{tabular}{|c|c|c|}
\multicolumn{3}{c}{\bf Table 1. Numerical Results: Small-Scale Finite} \\ 
\multicolumn{3}{c}{\bf Buffer Model With Exponential Arrivals} \\
\multicolumn{3}{c}{} \\
 \hline 
{\em p(n)} & Exponential-Direct  & Exponential-Traditional  \\       
\hline   
0 & 0.067958810 & 0.067958810 \\
1 & 0.169897026 & 0.169897026 \\
2 & 0.212371283 & 0.212371283 \\
3 & 0.176976069 & 0.176976069 \\
4 & 0.147480057 & 0.147480057 \\
5 & 0.122900048 & 0.122900048 \\
6 & 0.102416707 & 0.102416707 \\
\hline
\end{tabular}
\end{center}
\vspace{5mm}

%%%%% 

\begin{center}
\noindent
\begin{tabular}{|c|c|c|c|}
\multicolumn{4}{c}{\bf Table 2. Numerical Results: Small-Scale  } \\ 
\multicolumn{4}{c}{\bf Finite Buffer Model With Deterministic, } \\
\multicolumn{4}{c}{\bf  Erlang, or Hyper-Exponential Arrivals} \\
\multicolumn{4}{c}{} \\
 \hline 
{\em p(n)} & Deterministic & Erlang & Hyper-Exponential  \\       
\hline   
0 & 0.047234853 & 0.060394802 & 0.095667547 \\
1 & 0.127764093 & 0.153533580 & 0.170777140 \\
2 & 0.254669333 & 0.228194616 & 0.182173364 \\
3 & 0.232414603 & 0.196990341 & 0.153721590 \\
4 & 0.156638062 & 0.150739248 & 0.148862427 \\
5 & 0.107500964 & 0.117912665 & 0.131909935 \\
6 & 0.073778091 & 0.092234748 & 0.116887997 \\
\hline
\end{tabular}
\end{center}
\vspace{5mm}

%%%% Table 3 
\begin{center}
\noindent
\begin{tabular}{|c|c|c|}
\multicolumn{3}{c}{\bf Table 3. Performance Measures, Small-Scale Finite Buffer Models} \\ 
\multicolumn{3}{c}{} \\
 \hline 
 & $L$  & $W$  \\       
\hline   
Exponential - Direct       & 2.944488506 & 0.656092538 \\
Exponential - Traditional       & 2.944488506 & 0.656092538 \\
Deterministic - Direct     & 2.941072182 & 0.635068584 \\
Erlang - Direct            & 2.946822639 & 0.649247728 \\
Hyper-exponential - Direct  & 2.952616004 & 0.668684379 \\
\hline
\end{tabular}
\end{center}
\vspace{5mm}

\subsection{Large Size Applications: Computational Methodology and Results}

For the four distributions,  deterministic, Erlang,  exponential, and hyper-exponential,  we recursively compute the stationary distributions $\{\pi(.)\}$ and $\{p(.)\}$  using  the direct method as in Theorem~\ref{thm:pij} and the  algorithm in the Appendix.
 For the large scale examples, reported in Table 4 and Figures 2 and 3, we use $\lambda=5.8$, $\mu=.2$,  $c=30$, and $\rho=.96\overline{6}$. In the case of the hyper-exponential we use $p=.873563218, \lambda_1=8$, and $\lambda_2=2$ which gives a coefficient of variation $\approx 1.430035$ .  For each of the four distribution functions,  we use $\epsilon=10^{-125}$ in step 3 of the algorithm to recursively solve for $\sigma$, and use  Lemma~\ref{lem:n} (step  4 of the algorithm)  with $\epsilon= 10^{-16}$ to identify the truncation value  for the infinite capacity. The values for  $N$ are $490$, $705$, $918$ and $1349$ for the   deterministic, Erlang,  exponential, and hyper-exponential respectively.

We use the Python programming language to compute the stationary probabilities for large-scale examples. Use of the Decimal package, a fixed-decimal package capable of arbitrarily long mantissas, is notable in addressing the overflow errors associated with large factorials and the underflow issues created by the $LST$ values with large $c$ and $N$. 

The results for exponential inter-arrivals are compared with traditional methods for computing stationary probabilities using $M/M/c$ queues such as described in Gross \& Harris\cite{Gro08} and 
Kleinrock\cite{Kle75}.  This serves to verify, numerically, our method computations. % \par 
For deterministic, Erlang, and hyper-exponential inter-arrivals, our results are compared with Tak\'{a}cs method as generated by the QTS software provided by Gross \& Harris \cite{Gro08}.

As can be seen from the QTS (Tak\'{a}cs) results for the Erlang, deterministic, and hyper-exponential distributions, difficulties with floating-point overflow/underflow exist with this number of servers, and can persist as low as $c=10$.  These problems expand with increasing $c$, limiting usable results from that software.
\par
With the exception of $p(0)$, which compounds the error present in all other values of $p(n)$, our method is more numerically stable than Tak\'{a}cs as implemented by QTS software provided by Gross \& Harris \cite{Gro08} even when using floating-point levels of precision.  Our  use of the Decimal package provides accuracy for substantially higher values of $c$ while also reducing the error of $p(0)$ below that of floating-point implementations. 

Table 4 provides performance measures comparing our direct algorithm and Tak\'{a}cs method.  There is a good match between the two methods. Moreover, Figure 2  and Figure 3 provide the cumulative and density distribution functions  for  the four distributions, however we truncated the distribution functions in Figures 2 and 3 at $100$ and  $80$ respectively.

%%%% Table 3 
\begin{center}
\noindent
\begin{tabular}{|c|c|c|}
\multicolumn{3}{c}{\bf Table 4. Performance Measures} \\ 
\multicolumn{3}{c}{} \\
 \hline 
Distribution/Method & $L$  & $W$  \\       
\hline  
%%%

%%%% 
Deterministic - Direct    & 39.357&6.786\\
Deterministic - Tak\'{a}cs         & 39.357&6.786 \\
Erlang - Direct           & 45.600&7.870 \\
Erlang - Tak\'{a}cs                & 45.637&7.869 \\
Exponential - Direct      & 52.083&8.980 \\
Exponential - Traditional      & 52.083&8.980 \\
Hyper-exponential - Direct & 65.500&11.300 \\
Hyper-exponential - Tak\'{a}cs      & 65.480&11.290 \\
\hline
\end{tabular}
\end{center}

\begin{figure}[ht] % float placement: (h)ere, page (t)op, page (b)ottom, other (p)age  
 \includegraphics[bb=15 10 600 400,width=6.00in,height=4.00in,keepaspectratio]{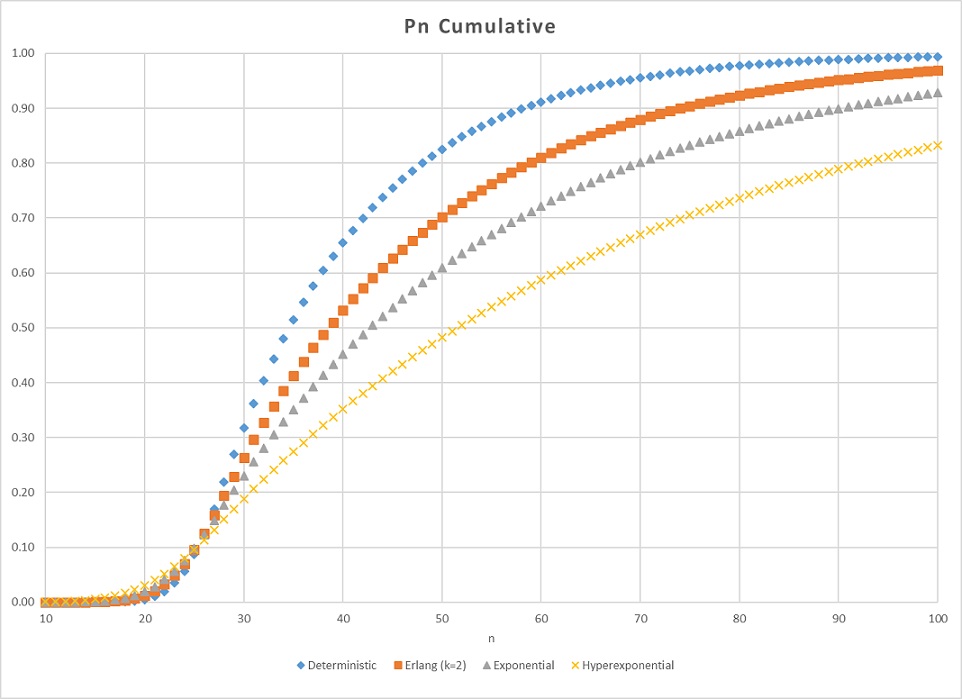}
  \caption{Cumulative Distribution Functions for $c=30, \rho=0.96\overline{6}$}
 \label{FigA1}
\end{figure}

\begin{figure}[ht] 
 \includegraphics[bb=15 10 600 400,width=6.00in,height=4.00in,keepaspectratio]{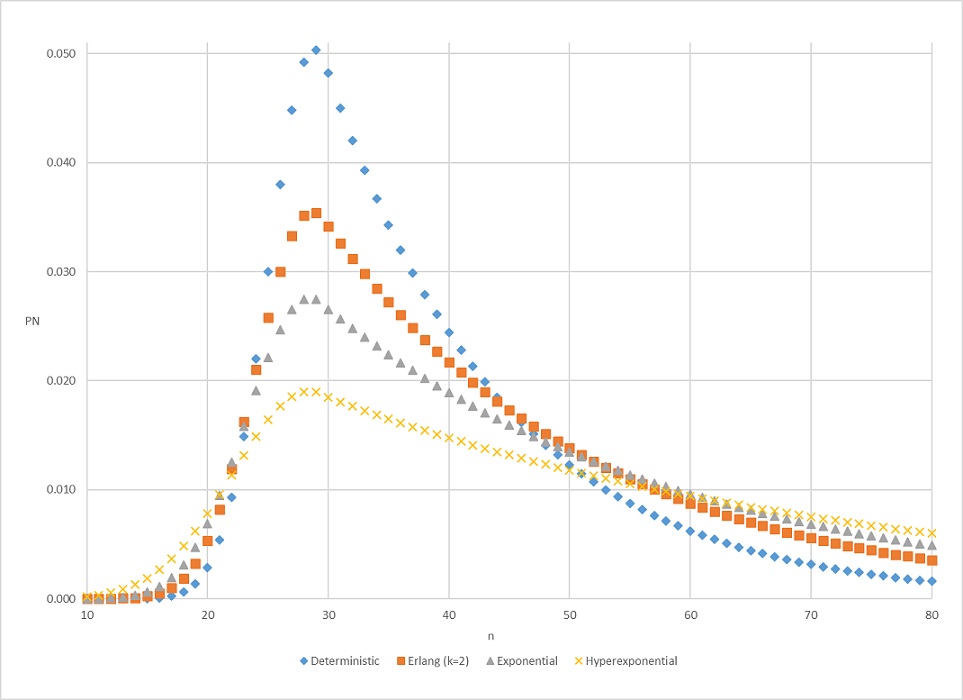}
  \caption{ Probability Mass Functions for $c=30, \rho=0.96\overline{6}$}
 \label{FigA2}
\end{figure}

\bigskip

We do not give numerical results for the batch arrival models to save space. Our approach in dealing with batch arrival model is to build on the corresponding  standard $G/M/c$ model as indicated in the algorithm given in Subsection  ~\ref{subsec:adapt} of the appendix.

\pagebreak
.
\newpage
.

\section{Appendix A}
The direct algorithm presented here is used to compute the stationary probabilities and performance measures for both the finite (small scale) and the infinite (large scale) buffer cases.
Our computational methodology follows that provided in Sections 2 and 3, wherein the transition probabilities are  prepared following the examples described in the first part of Section 4. Moreover, $|\sigma|<1$ is determined from $A^*(c\mu(1-\sigma))=\sigma$, and $\pi(j)$ values are prepared using $\sigma$ and the transition probabilities, and then normalized to provide steady-state system size probabilities from which we obtain the performance measures. Below we give an algorithm  for the multi-server model, followed by an adaptation to the batch-arrival multi-server case. 

%\noindent
%{\bf Algorithm.}
\subsection{Algorithm for the $GI/M/c$ Model}
 Initialization. Let $\epsilon$ be  the maximum allowable error,  $N$ be the finite buffer for small examples,
 $c$ be the number of servers each with mean rate $\mu$, and $\lambda$ be the mean arrival rate. We also input the $LST$ of the inter-arrival times distribution and number of phases and/or weights if applicable. Note that for large scale examples we determine $N$ using  $\epsilon$. \par
%%%%%%%
\ls{1}
\begin{enumerate}
 %%%%%%%%%%%%%%
	\item Compute $\rho = \lambda/ c\mu$ and check that $\rho <1$ (i.e.: a long-run solution exists).        
	\item For each specified inter-arrival distribution, compute\; $A^*(s)$ for $s=k\mu$ where $k=1,2, ..., c$ and $A^*_n(c\mu)$ for $n=1,2, ..., N-c+1$\\\\
  Deterministic: 
$A^*(s) = e^{-s/\lambda}$,\;\;    $  A^*_n(c\mu) = \lambda^{-n} e^{-c\mu/ \lambda}$ \\ \\
	 Exponential: 
	$A^*(s) = \frac{\lambda}{s+\lambda} $,\;\; $ A^*_n(c\mu) = \frac{n! \lambda}{(c\mu+\lambda)^{n+1}}$ \\\\
	Erlang (two-phase):
	 $A^*(s) = \frac{4\lambda^2}{(s+2\lambda)^{2}} $,\;\; $  A^*_n(c\mu) = n!(n+1) \frac{4\lambda^2}{(c\mu+2\lambda)^{n+2}} $\\\\
Hyper-exponential (two phase): 
	 $A^*_n(s)= n! \left[ \frac{p\lambda_1}{(s+\lambda_1)^{n+1}} + \frac{(1-p)\lambda_2}{(s+\lambda_2)^{n+1}}\right]\;, \text{ and}  $\\ 
        $ A^*_n(c\mu)= n! \left[ \frac{p\lambda_1}{(c\mu+\lambda_1)^{n+1}} + \frac{(1-p)\lambda_2}{(c\mu+\lambda_2)^{n+1}}\right]\;. $
%
%%%%%%%%%%%%%%%%%%%%%%%
	\item Root-solve by iterating over the following until 
	$| \sigma_{n+1}-\sigma_{n}| < \epsilon$\\
	\[\sigma_{n+1}=A^*[c\mu (1-\sigma_n)]\]
\item For a given $\epsilon$, determine $N$ (for large examples) using
	\[
 N = \min \left\{ n \in {\cal N} \; \big | \; n \geq c+\frac{\ln(\epsilon)-2\ln(1-\sigma)}{\ln(\sigma)} \right\}
 \]	
%%%%%%%%%%%%%%%%
	\item Define $p(i,j) = 0$ for all $i = 0, 1,...,N-2$,\;\, $j = i+2, i+3, ..., N$.
	\item Compute $p(i,j)$ for $i=0,1,2,...,c-1$,\;\, $j=1,2,...,i+1$ using the computationally friendly  form:
	\[p(i,j) = \frac{(i+1)!}{j!} \sum_{k=0}^{i-j+1} \frac{(-1)^k A^*((j+k)\mu)}{(i-j-k+1)!k!}\;. \]
	\item Compute $p(i,j)$ for $i=c,c+1,...,N$,\;\, $j=c,c+1,...,i+1, \;\;\; i+1 \leq N$ using
	\[
	p(i,j)=  \frac {(c\mu )^{i-j+1}A^*_{i-j+1}(c\mu)}{(i-j+1)!}
	\]
	\item Compute $p(i,j)$ for $i=c,c+1,...,N$,\;\, $j=1,2,...,c-1$ using the computationally friendly  form:
\begin{align*}
p(i,j) 
&= \frac{c!}{j!}\sum^{c-j}_{k=1} \frac{(-1)^{(c-j-k)}}{k!(c-j-k)!} \left(\frac{c}{k}\right)^{i-c+1} 
\left[ A^*((c-k)\mu)- \sum^{i-c+1}_{r=0} \frac{(k\mu )^{r}A^*_r(c\mu)}{r!}  \right] \; 
\end{align*}
	\item Compute $p(i,j)$ for $i=c,c+1,...,N$,\;\, $j=0$ using
	\[p(i,j) = 1 - \sum_{n=1}^{N}p_{i,n} \]
	\item Define $a(k,j) = 0$ for $k=0,1,..., N$,\;\, $j = k, k+1, ..., N$  
	\item Compute $a(k,j)$ for $k=j+1,j+2,..., N$,\;\, $j=0,1,...,c-1$ using
	\[a(k,j)=\sum_{i=0}^{j}p(k,i)\]
	\item Compute $\pi' (j) = \sigma^j$ for $j=c, c+1, ..., N$
	\item Compute $\pi' (j) $ for $j=c-1,c-2,...,0$ recursively using
	\[\pi'(j) = \frac{\sum_{k=j+1}^N\pi'(k)a(k,j)}{p(j,j+1)}\]
	\item Compute $\pi(j)$ for $j=0,1,...,N$ by normalizing $\pi'_j$ using
	\[\pi(j)=\frac{\pi'(j)}{\Phi} \] where \[\Phi= {\sum_{k=0}^N\pi'(k)} = {\sum_{k=0}^{c-1}\pi'(k) + \sum_{k=c}^{N}\sigma^k} = {\sum_{k=0}^{c-1}\pi'(k) + \frac{\sigma^c(1-\sigma^{N-c+1})}{(1-\sigma)}}\]
	\item Compute $p(0)$ using
	\[ p(0) = (1-\rho) + \rho\pi(N) - \rho\sum_{k=0}^{c-2}\frac{c-k-1}{k+1}\pi(k) \]
	\item Compute $p(n)$ for $n = 1, 2, ..., c-1$ using
	\[p(n) = \frac{c\rho \pi(n-1)}{n}\]
	\item Compute $p(n)$ for $n = c, c+1, ..., N$ using
	\[p(n) = \rho \pi(n-1)\]
	\item Compute performance measures using
	$E[L]=\sum_{i=1}^{N}ip(i)$ and $ E[W]=E[L]/\lambda$.	
\end{enumerate}
\ls{1.3}

\subsection{Adaptation to the $GI^X/M/c/N$ model}\label{subsec:adapt}

 Here, we identify steps in Appendix A  that need to be modified.
We refer to cases 2 and 3 in Section~\ref{sec:extensions}.
\begin{enumerate}
\item Repeat steps 1 and 2 in Appendix A
\item  Compute $p^*(i,j)$ as in Section~\ref{sec:extensions} for case 2 and case 3
\item Compute arrival-time probabilities $\{\pi(i),i=0,1,\ldots ,N\}$ by solving  $\pi=\pi P^*$, $\sum_{i=0}^N \pi(i)=1$. 

\item  Compute $\{p(i)\}'s$  for cases 2 and 3 using  Section~\ref{sec:extensions}. 
\end{enumerate}

%\section{Declarations}
%The authors declare that they have no conflict of interest.
%
%This article received no external funding.
%
%All data used during the current study is generated by Python code. The code is available upon request.
%
%%\bibliography{ref}

\begin{thebibliography}{10}

\bibitem{Bai81}
D.~E. Baily and M.~F. Neuts.
\newblock Algorithmic methods for multi-server queues with group arrivals and
  exponential services.
\newblock {\em European Journal of Operational Research}, 8(2):184--196, 1981.

\bibitem{Cha16}
M.L. Chaudhry and J.J. Kim.
\newblock Analytically elegant and computationally efficient results in terms
  of roots for the {$GI^X/M/c$} queueing system.
\newblock {\em Queueing Systems}, 82:237--257, 2016.

\bibitem{Cos80}
G.~Cosmetatos and S.~Godsave.
\newblock Approximations in the multi-server queue with hyper-exponential
  inter-arrival times and exponential service times.
\newblock {\em Journal of the Operational Research Society}, 31(1):57--62,
  1980.

\bibitem{Elt21}
M.~El-Taha.
\newblock An efficient convolution method to compute the stationary transition
  probabilities of the g/m/c model and its variants.
\newblock In {\em 2021 10th IFIP International Conference on Performance
  Evaluation and Modeling in Wireless and Wired Networks (PEMWN)}, pages 1--6.
  IEEE, 2021.

\bibitem{Elt99}
M.~El-Taha and S.~Stidham~Jr.
\newblock {\em Sample-Path Analysis of Queueing Systems}.
\newblock Kluwer Academic Publishing, Boston, 1999.

\bibitem{Fer06}
F.~Ferreira and A.~Pacheco.
\newblock Analysis of {GI}/{M}/s/c queues using uniformisation.
\newblock {\em Computers and Mathematics with Applications}, 51:291--304, 2006.

\bibitem{Gra14}
W.~Grassmann and J.~Tavakoli.
\newblock Efficient methods to find the equilibrium distribution of the number
  of customers in {GI/M/c} queues.
\newblock {\em IFOR}, 52:197--205, 2014.

\bibitem{Gro08}
D.~Gross, J.F. Shortle, J.M. Thompson, and C.~Harris.
\newblock {\em Fundamentals of Queueing Theory}.
\newblock John Wiley, New Jersey, 4th edition, 2008.

\bibitem{Hok75}
Per Hokstad.
\newblock The {G/M/m} queue with finite waiting room.
\newblock {\em Journal of Applied Probability}, 12(4):779--792, 1975.

\bibitem{Kel79}
F.~Kelly.
\newblock {\em Reversibility and Stochastic Networks}.
\newblock Wiley, New York, 1979.

\bibitem{Kim17}
J.J. Kim and M.L. Chaudhry.
\newblock A novel way of treating the finite-buffer queue ${GI/M/c/N}$ using
  roots.
\newblock {\em International Journal of Mathematical Models and Methods in
  Applied Sciences}, 11, 2017.

\bibitem{Kle75}
L.~Kleinrock.
\newblock {\em Queueing {S}ystems, vol. I}.
\newblock Wiley Intersciences, New York, 1975.

\bibitem{Lax00}
P.~V. Laxmi and U.C. Gupta.
\newblock Analysis of finite-buffer multi-server queues with group arrivals:
  {$GI^X/M/c/N$}.
\newblock {\em Queueing Systems}, 36(1-3):125--140, 2000.

\bibitem{Med03}
J.~Medhi.
\newblock {\em Stochastic Models in Queueing Theory}.
\newblock Academic Press, New York, 2nd edition, 2003.

\bibitem{Neu81}
M.F. Neuts.
\newblock {\em Matrix-Geometric Solutions in Stochastic Models-An Algorithmic
  Approach}.
\newblock Johns Hopkins University Press, 1981.

\bibitem{Ros07}
S.M. Ross.
\newblock {\em Introduction to Probability Models}.
\newblock Academic Press, San Diego, 9th edition, 2007.

\bibitem{Sti87}
S.~Stidham~Jr.
\newblock Stable recursive procedures for numerical computations in {M}arkov
  models.
\newblock {\em Annals of Operations Research}, 8:27--40, 1987.

\bibitem{Sti89}
S.~Stidham~Jr. and M.~El-Taha.
\newblock Sample-path analysis of processes with imbedded point processes.
\newblock {\em Queueing Systems}, 5:131--165, 1989.

\bibitem{Tak57}
L.~Tak\'{a}cs.
\newblock On a queueing problem concerning telephone traffic.
\newblock {\em Acta Math.Acad.Sci. Hungar.}, 8:325--335, 1957.

\bibitem{Tak62}
L.~Tak\'{a}cs.
\newblock {\em Introduction to the theory of queues}.
\newblock Oxford University Press, New York, 1962.

\bibitem{Tij03}
H.C. Tijms.
\newblock {\em A First Course in Stochastic Models}.
\newblock Wiley, New York, 2003.

\bibitem{Yao84}
D.D. Yao, M.L. Chaudhry, and J.G.C. Templeton.
\newblock A note on some relations in the queue {$GI^X/M/c$}.
\newblock {\em Operations Research Letters}, 3(1):53--56, 1984.

\end{thebibliography}

%\end{document}

\end{document}